\documentclass[12pt]{article}
\usepackage{amssymb,amsmath}
\usepackage{cases}
\usepackage{amsfonts}
\usepackage{color,xcolor,cite}
\usepackage[left=2.0cm,right=2.0cm,top=2.0cm,bottom=2.0cm]{geometry}
\usepackage[colorlinks,citecolor=blue,urlcolor=blue]{hyperref}

\newtheorem{theorem}{Theorem}[section]

\newtheorem{lemma}{Lemma}[section]
\newtheorem{proposition}{Proposition}[section]

\newtheorem{definition}{Definition}[section]
\newtheorem{remark}{Remark}[section]

\newcommand{\bal}{\begin{align}}
\newcommand{\bbal}{\begin{align*}}
\newcommand{\beq}{\begin{equation}}
\newcommand{\eeq}{\end{equation}}
\newcommand{\bca}{\begin{cases}}
\newcommand{\eca}{\end{cases}}

\newcommand{\pa}{\partial}
\newcommand{\fr}{\frac}

\newcommand{\De}{\Delta}

\newcommand{\dd}{\mathrm{d}}

\newcommand{\R}{\mathbb{R}}
\newcommand{\T}{\mathbb{T}}
\newcommand{\Z}{\mathbb{Z}}
\newcommand{\vv}{\mathbf{v}}

\newcommand{\bi}{\Big}
\newcommand{\g}{\big}
\begin{document}
\title{Non-uniform dependence on initial data for the Camassa-Holm equation in the critical Besov space}

\author{Jinlu Li$^{1}$\footnote{E-mail: lijinlu@gnnu.edu.cn}, Xing Wu$^2$\footnote{E-mail: ny2008wx@163.com}, Yanghai Yu$^{3}$\footnote{E-mail:  yuyanghai214@sina.com(Corresponding author)} and Weipeng Zhu$^{4}$\footnote{E-mail: mathzwp2010@163.com}\\
\small $^1$ School of Mathematics and Computer Sciences, Gannan Normal University, Ganzhou 341000, China\\
\small $^2$ College of Information and Management Science, Henan Agricultural University, Zhengzhou 450002, China\\
\small $^3$ School of Mathematics and Statistics, Anhui Normal University, Wuhu 241002, China\\
\small $^4$ School of Mathematics and Information Science, Guangzhou University, Guangzhou 510006, China}

\date{\today}

\maketitle\noindent{\hrulefill}

{\bf Abstract:} Whether or not the data-to-solution map of the Cauchy problem for the Camassa-Holm equation and Novikov equation in the critical Besov space $B_{2,1}^{3/2}(\R)$ is not uniformly continuous remains open. In the paper, we aim at solving the open question left the previous works in \cite{Li3,Li4} and give a positive answer to this problem.

{\bf Keywords:} Camassa-Holm (Novikov) equation, Non-uniform continuous dependence, Critical Besov spaces

{\bf MSC (2010):} 35Q35; 35A01; 76W05
\vskip0mm\noindent{\hrulefill}

\section{Introduction}

In this paper, we consider the Cauchy problem for the well-known Camassa-Holm equation
\begin{equation*}
\rm{(CH)}\quad\begin{cases}
u_t-u_{xxt}+3uu_x=2u_xu_{xx}+uu_{xxx}, \; &(x,t)\in \R\times\R^+,\\
u(x,t=0)=u_0,\; &x\in \R.
\end{cases}
\end{equation*}

(CH) was firstly proposed in the context of hereditary symmetries studied in \cite{Fokas} and then was derived explicitly as a water wave equation by Camassa--Holm \cite{Camassa}. (CH) is completely integrable \cite{Camassa,Constantin-P} with a bi-Hamiltonian structure \cite{Constantin-E,Fokas} and infinitely many conservation laws \cite{Camassa,Fokas}. Also, it admits exact peaked
soliton solutions (peakons) of the form $ce^{-|x-ct|}$ with $c>0$, which are orbitally stable \cite{Constantin.Strauss} and models wave breaking (i.e., the solution remains bounded, while its slope becomes unbounded in finite
time \cite{Constantin,Escher2,Escher3}. It is worth mentioning that the peaked solitons present the characteristic for the travelling water waves of greatest height and largest amplitude and arise as solutions to the free-boundary problem for incompressible Euler equations over a flat bed, see Refs. \cite{Constantin-I,Escher4,Escher5,Toland} for the details. Because of the mentioned interesting and remarkable features,
the CH equation has attracted much attention as a class of integrable shallow water wave equations in recent twenty years. Its systematic mathematical study was initiated in a series of papers by Constantin and Escher, see \cite{Escher1,Escher2,Escher3,Escher4,Escher5}. We can refer the readers to see the global strong solutions in \cite{Constantin,Escher1,Escher2} and finite time blow-up strong solutions in \cite{Constantin,Escher1,Escher2,Escher3} to (CH), the existence and uniqueness of global weak solutions in \cite{Constantin.Molinet, Xin.Z.P}, the global conservative solutions \cite{Bre1} and global dissipative solutions \cite{Bre2} in $H^1(\R)$.

After the phenomenon of non-uniform continuity for some dispersive equations was studied by Kenig et al. \cite{Kenig2001}, the issue of non-uniform dependence on the initial data has been a fascinating object of research in the recent past. Naturally, we may wonder which regularity assumptions are relevant for the initial data $u_0$ such that the Cauchy problem to (CH) is not uniform dependence on initial data, namely, the dependence of solution on the initial data associated with this equation is not uniformly continuous. Himonas--Misio{\l}ek \cite{H-M} obtained the first result on the non-uniform dependence for (CH)  in $H^s(\T)$ with $s\geq2$ using explicitly constructed travelling wave solutions, which was sharpened to $s>\fr32$ by Himonas--Kenig \cite{H-K} on the real-line and Himonas--Kenig--Misio{\l}ek \cite{H-K-M} on the circle. Danchin \cite{d1,d3} proved the local existence and uniqueness of strong solutions to (CH) with initial data in $B^s_{p,r}$ if $(p,r)\in[1,\infty]\times[1,\infty), s>\max\g\{1+\fr1p, \fr32\g\}$ and $B^{3/2}_{2,1}$. Li--Yin \cite{Li-Yin1} proved that the continuity of the solution map of (CH) with respect to the initial data.  Guo et al.\cite{Guo-Yin} established the ill-posedness of (CH) in $H^{3/2}$ and in $B_{2,r}^{3/2}$ with $r\in(1,\infty)$ by proving the norm inflation. In our recent paper\cite{Li3}, we proved that the non-uniform dependence on initial data for (CH) under the framework of Besov spaces $B^s_{p,r}$ for $s>\max\g\{1+\fr1p, \fr32\g\}$. However, whether or not the data-to-solution map of the Cauchy problem for (CH) in the critical Besov spaces $B_{2,1}^{3/2}(\R)$ is not uniformly continuous remains open. We aim at giving a positive answer to this question in this paper.

Before stating our main result, we transform (CH) equivalently into the following nonlinear transport type equation
\begin{equation}\label{CH}
\begin{cases}
\partial_tu+u\pa_xu=\mathbf{P}(u), \; &(x,t)\in \R\times\R^+,\\
u(x,t=0)=u_0,\; &x\in \R,
\end{cases}
\end{equation}
where
\begin{equation}\label{CH1}
\mathbf{P}(u)=P(D)\bi(u^2+\fr12(\pa_xu)^2\bi)\quad\text{with}\quad P(D)=-\pa_x\g(1-\pa^2_x\g)^{-1}.
\end{equation}
Our main result is stated as follows.
\begin{theorem}\label{the1.1}
The solution map $u_0\rightarrow \mathbf{S}_t(u_0)$ of the Cauchy problem \eqref{CH}--\eqref{CH1} is not uniformly continuous from any bounded subset in $B^{\frac32}_{2,1}$ into $\mathcal{C}([0,T];B^{\frac32}_{2,1})$. More precisely, there exists two sequences of solutions $\mathbf{S}_t(f_n+g_n)$ and $\mathbf{S}_t(f_n)$ such that
\bbal
&||f_n||_{B^{\frac32}_{2,1}}\lesssim 1 \quad\text{and}\quad \lim_{n\rightarrow \infty}||g_n||_{B^{\frac32}_{2,1}}= 0
\end{align*}
but
\bbal
\liminf_{n\rightarrow \infty}||\mathbf{S}_t(f_n+g_n)-\mathbf{S}_t(f_n)||_{B^{\frac32}_{2,1}}\gtrsim t,  \quad \forall \;t\in[0,T_0],
\end{align*}
with small time $T_0$.
\end{theorem}
\begin{remark}\label{re0}
The method we used in \cite{Li3} does not work for the critical index $s=\fr32$ due to technical difficulty which mainly lies in the transport equation theory forbids the estimate of solution in $B_{2,1}^{1/2}$. We prove Theorem \ref{the1.1} by utilizing new method and generalize the previous result \cite{Li3} to the critical case.
\end{remark}
\begin{remark}\label{re1}
The method we used in proving the Theorem \ref{the1.1} can be applied equally well to other related systems,
such as the following Novikov equation
\begin{equation}\label{novikov}
\begin{cases}
u_t+u^2u_x=\mathbf{Q}(u),\\
u(x,t=0)=u_0,
\end{cases}
\end{equation}
where
\begin{equation}\label{novikov1}
\mathbf{Q}(u)=-(1-\pa^2_x)^{-1}\Big(\frac12u_x^3+\pa_x\big(\frac32uu^2_x+u^3\big)\Big).
\end{equation}
\end{remark}
Then we have the following
\begin{theorem}\label{the1.2}
The solution map $u_0\rightarrow \mathbf{S}_t(u_0)$ of the Cauchy problem \eqref{novikov}--\eqref{novikov1} is not uniformly continuous from any bounded subset in $B^{\frac32}_{2,1}$ into $\mathcal{C}([0,T];B^{\frac32}_{2,1})$. More precisely, there exists two sequences of solutions $\mathbf{S}_t(f_n+h_n)$ and $\mathbf{S}_t(f_n)$ such that
\bbal
&||f_n||_{B^{\frac32}_{2,1}}\lesssim 1 \quad\text{and}\quad \lim_{n\rightarrow \infty}||h_n||_{B^{\frac32}_{2,1}}= 0
\end{align*}
but
\bbal
\liminf_{n\rightarrow \infty}||\mathbf{S}_t(f_n+h_n)-\mathbf{S}_t(f_n)||_{B^{\frac32}_{2,1}}\gtrsim t,  \quad \forall \;t\in[0,T_0],
\end{align*}
with small time $T_0$.
\end{theorem}
\noindent\textbf{Organization of our paper.} In Section 2, we list some notations and known results which will be used in the sequel. In Section 3, we present the local well-posedness result and establish some technical Propositions. In Section 4, we prove our main theorem by adopting the strategies used in \cite{Li3}. Here we should point out the new difficulty when dealing the critical case lies in the lack of the estimate of solution in $B^{\frac12}_{2,1}$. To overcome this, we decompose the solution map as
 \bbal
\mathbf{S}_{t}(u^n_0)=\mathbf{S}_{t}(u^n_0)-u^n_0-t\vv_0(u_0^n)+f_n+g_n+t\big(\mathbf{P}(u^n_0)-u^n_{0}\pa_xu^n_{0}\big).
\end{align*}
On one hand, $u^n_{0}\pa_xu^n_{0}$ brings us the term $g_n\partial_xf_n$ which plays an essential role since it would not small when $n$ is large enough; On the other hand, $\mathbf{S}_{t}(u^n_0)-u^n_0-t\vv_0(u_0^n)$ promotes us to estimate the crucial quantity $\|\mathbf{S}_{t}(u^n_0)\|_{L^\infty}$ which can be controlled by $t\|u^n_0\|^2_{C^{0,1}}+\|u^n_0\|_{L^\infty}$ instead of the norm $\|\mathbf{S}_{t}(u^n_0)\|_{B^{\frac12}_{2,1}}$.
Based on the suitable choice of $f_n$ and $g_n$, we prove that the solution map is not uniformly continuous.

\section{Littlewood-Paley analysis}
We firstly introduce some notations which will be used throughout this paper.

The symbol $A\lesssim(\gtrsim) B$ means that there is a uniform positive constant $c$ independent of $A$ and $B$ such that $A\leq(\geq) cB$. Given a Banach space $X$, we denote its norm by $\|\cdot\|_{X}$. We use the simplified notation $||\mathbf{f}_1,\cdots,\mathbf{f}_n||_X=||\mathbf{f}_1||_X+\cdots+||\mathbf{f}_n||_X$ if without confusion.
For all $f\in \mathcal{S}'$, the Fourier transform $\mathcal{F}f$ (also denoted by $\hat{f}$) is defined by
$$
\mathcal{F}f(\xi)=\hat{f}(\xi)=\int_{\R}e^{-ix\xi}f(x)\dd x \quad\text{for any}\; \xi\in\R.
$$
The inverse Fourier transform allows us to recover $u$ from $\hat{f}$:
$$
f(x)=\mathcal{F}^{-1}\hat{f}(x)=\frac{1}{2\pi}\int_{\R}e^{ix\xi}\hat{f}(\xi)\dd\xi.
$$
Next, we will recall some facts about the Littlewood-Paley decomposition, the nonhomogeneous Besov spaces and their some useful properties (see \cite{B.C.D} for more details).

There exists a couple of smooth functions $(\chi,\varphi)$ valued in $[0,1]$, such that $\chi$ is supported in the ball $\mathcal{B}\triangleq \{\xi\in\mathbb{R}:|\xi|\leq \frac 4 3\}$, and $\varphi$ is supported in the ring $\mathcal{C}\triangleq \{\xi\in\mathbb{R}:\frac 3 4\leq|\xi|\leq \frac 8 3\}$. Moreover,
\begin{eqnarray*}
\chi(\xi)+\sum_{j\geq0}\varphi(2^{-j}\xi)=1 \quad \mbox{ for any } \xi\in \R.
\end{eqnarray*}

For every $f\in \mathcal{S'}(\mathbb{R})$, the inhomogeneous dyadic blocks ${\Delta}_j$ are defined as follows
\begin{numcases}{\Delta_jf=}
0, & if $j\leq-2$;\nonumber\\
\chi(D)f=\mathcal{F}^{-1}(\chi \mathcal{F}f), & if $j=-1$;\nonumber\\
\varphi(2^{-j}D)f=\mathcal{F}^{-1}\g(\varphi(2^{-j}\cdot)\mathcal{F}f\g), & if $j\geq0$.\nonumber
\end{numcases}
In the inhomogeneous case, the following Littlewood-Paley decomposition makes sense
$$
f=\sum_{j\geq-1}{\Delta}_jf\quad \text{for any}\;f\in \mathcal{S'}(\mathbb{R}).
$$

\begin{definition}[See \cite{B.C.D}]
Let $s\in\mathbb{R}$ and $(p,r)\in[1, \infty]^2$. The nonhomogeneous Besov space $B^{s}_{p,r}(\R)$ is defined by
\begin{align*}
B^{s}_{p,r}(\R):=\Big\{f\in \mathcal{S}'(\R):\;\big\|\g(2^{js}\|\Delta_j{f}\|_{L^p(\R)}\g)_{j\in \Z}\big\|_{\ell^r(\Z)}<\infty\Big\}.
\end{align*}
\end{definition}
Finally, we give some important properties which will be also often used throughout the paper.
\begin{lemma}[See \cite{B.C.D}]\label{le2}
For $s>0$, then for any $u,v \in B^{s}_{2,1}(\R)\cap L^\infty(\R)$, we have
\bbal
&\|uv\|_{B^{s}_{2,1}(\R)}\leq C\big(\|u\|_{B^{s}_{2,1}(\R)}\|v\|_{L^\infty(\R)}+\|v\|_{B^{s}_{2,1}(\R)}\|u\|_{L^\infty(\R)}\big).
\end{align*}
In particular, we have the embedding $B^{1/2}_{2,1}(\R)\hookrightarrow L^\infty(\R)$ and
$$B^s_{p,q}(\R)\hookrightarrow B^t_{p,r}(\R)\quad\text{for}\;s>t\quad\text{or}\quad s=t,1\leq q\leq r\leq\infty.$$
\end{lemma}

\begin{lemma}[Lemma 3.26 in \cite{B.C.D}]\label{01}
Let $(p,r)\in[1, \infty]^2$, $s>1$ and $u_0\in B^s_{p,r}(\R)$.
Assume that $u\in L^\infty([0,T]; B^s_{p,r}(\R))$ solves \eqref{CH}--\eqref{CH1}.
Then there exists a constant $C=C(s,p)$ and a universal constant $C'$ such that for all $t\in[0,T]$, we have
\begin{align*}
&||u(t)||_{B^s_{p,r}(\R)}\leq ||u_0||_{B^s_{p,r}(\R)}\exp\Big\{C\int_0^t \|u(\tau)||_{C^{0,1}(\R)}\mathrm{d}\tau\Big\},\\
&||u(t)||_{C^{0,1}(\R)}\leq ||u_0||_{C^{0,1}(\R)}\exp\Big\{C'\int_0^t \|\pa_xu(\tau)||_{L^{\infty}(\R)}\mathrm{d}\tau\Big\}.
\end{align*}
\end{lemma}

Let us recall the local well-posedness result for (CH) in the critical Besov spaces.

\begin{lemma}[See \cite{d3}]\label{le5}
For any initial data $u_0$ which belongs to $$B_R=\g\{\psi\in B_{2,1}^\fr32: ||\psi||_{B^{\fr32}_{2,1}}\leq R\g\}\quad\text{for any}\;R>0.$$ Then there exists some $T=T\big(||u_0||_{B_{2,1}^\fr32}\big)>0$ such that (CH) has a unique solution $\mathbf{S}_{t}(u_0)\in \mathcal{C}([0,T];B^{\fr32}_{2,1})$. Moreover, we have
\begin{align*}
||\mathbf{S}_{t}(u_0)||_{B^\fr32_{2,1}}\leq C||u_0||_{B_{2,1}^\fr32}.
\end{align*}
\end{lemma}

\section{The Key Estimations}
Firstly, we need to introduce smooth, radial cut-off functions to localize the frequency region. Precisely,
let $\hat{\phi}\in \mathcal{C}^\infty_0(\mathbb{R})$ be an even, real-valued and non-negative function on $\R$ and satisfy
\begin{numcases}{\hat{\phi}(\xi)=}
1,&if $|\xi|\leq \frac{1}{4}$,\nonumber\\
0,&if $|\xi|\geq \frac{1}{2}$.\nonumber
\end{numcases}
Next, we need to establish the following crucial lemmas which will be used later on.
\begin{lemma}\label{le4} Let $(p,r)\in [1,\infty]\times[1,\infty)$. We define the high frequency function $f_n$ and the low frequency functions $g_n,h_n$ as follows
\bbal
&f_n=2^{-\fr32n}\phi(x)\sin \bi(\frac{17}{12}2^nx\bi),\\
&g_n=\frac{12}{17}2^{-n}\phi(x)\quad\text{and}\\
&h_n=\frac{12}{17}2^{-\frac{n}{2}}\phi(x),\quad n\gg1.
\end{align*}
Then for any $\sigma\in\R$, we have
\bal
&\|f_n\|_{L^\infty}\leq C2^{-\fr{3}{2}n}\phi(0)\quad\text{and}\quad||\pa_xf_n||_{L^\infty}\leq C2^{-\fr{n}{2}}\phi(0),\label{y00}\\
&\|g_n,\pa_xg_n\|_{L^\infty}\leq C2^{-n}\phi(0)\quad\text{and}\quad\|h_n,\pa_xh_n\|_{L^\infty}\leq C2^{-\fr{n}{2}}\phi(0),\label{y0}\\
&\|g_n\|_{B^\sigma_{p,r}}\leq C2^{-(n+\sigma)}\|\phi\|_{L^p}\quad\text{and}\quad \|h_n\|_{B^\sigma_{p,r}}\leq C2^{-(\frac{n}{2}+\sigma)}\|\phi\|_{L^p},\label{y1}\\
&\|f_n\|_{B^\sigma_{p,r}}\leq C2^{(\sigma-\fr32)n}\|\phi\|_{L^p},\label{y2}\\
&\liminf_{n\rightarrow \infty}\|g_n\pa_xf_n\|_{B^{\fr32}_{2,\infty}}\geq M_1,\label{y3}\\
&\liminf_{n\rightarrow \infty}\|h^2_n\pa_xf_n\|_{B^{\fr32}_{2,\infty}}\geq M_2,\label{y4}
\end{align}
for some positive constants $C, M_1,M_2$.
\end{lemma}
{\bf Proof.}\quad Direct computations gives \eqref{y00}--\eqref{y0}. Notice that
\bbal
\mathrm{supp} \ \hat{g}_n\subset \Big\{\xi\in\R: \ 0\leq |\xi|\leq \fr12\Big\},
\end{align*}
then, we have
\bbal
\widehat{\Delta_jg_n}=\varphi(2^{-j}\xi)\hat{g}_n(\xi)\equiv0\quad\text{for}\quad j\geq0,
\end{align*}
which implies
\bbal
{\Delta_jg_n}\equiv0\quad\text{for}\quad j\geq0.
\end{align*}
By the definitions of $g_n$ and the Besov space, we deduce that
\bbal
||g_n||_{B^\sigma_{p,r}}&=\frac{12}{17}2^{-(n+\sigma)}||\De_{-1}\phi||_{L^p}\leq \frac{12}{17}2^{-(n+\sigma)}||\phi||_{L^p}.
\end{align*}
Notice that
\bbal
\mathrm{supp} \ \widehat{h_n}\subset \Big\{\xi\in\R: \ 0\leq |\xi|\leq \fr12\Big\}\;\Rightarrow\;\mathrm{supp} \ \widehat{h^2_n}\subset \Big\{\xi\in\R: \ 0\leq |\xi|\leq 1\Big\},
\end{align*}
then, we have
\bbal
\mathrm{supp}\ \widehat{h^2_n\pa_xf_n}\subset \Big\{\xi\in\R: \ \frac{17}{12}2^n-\fr32\leq |\xi|\leq \frac{17}{12}2^n+\fr32\Big\},
\end{align*}
which implies
\begin{numcases}{\Delta_j\g(h^2_n\pa_xf_n\g)=}
h^2_n\pa_xf_n, &if $j=n$,\nonumber\\
0, &if $j\neq n$.\nonumber
\end{numcases}
By the definitions of $f_n$ and $h_n$, we obtain for some $\delta>0$
\bbal
||h^2_n\pa_xf_n||_{B^{\fr32}_{2,\infty}}&=2^{\fr32n}||\De_{n}\g(h^2_n\pa_xf_n\g)||_{L^2}=2^{\fr32n}||h^2_n\pa_xf_n||_{L^2}
\\&=\bi\|\frac{12}{17}\phi^3(x)\cos \bi(\frac{17}{12}2^nx\bi)+\Big(\frac{12}{17}\Big)^22^{-n}\phi^2(x)\pa_x\phi(x)\sin \bi(\frac{17}{12}2^nx\bi)\bi\|_{L^2}
\\&\geq \frac{12}{17}\bi\|\phi^3(x)\cos \bi(\frac{17}{12}2^nx\bi)\bi\|_{L^2}-C2^{-n}\\
&\geq \frac{1}{8}\cdot\frac{12}{17}\delta\phi^3(0)\bi(\int^\delta_{0}\bi|\cos\bi(\fr{17}{12}2^nx\bi)\bi|^2\dd x\bi)^{1/2}-C2^{-n}.
\end{align*}
We thus deduce that \eqref{y4}(see Lemma 1 in \cite{Li3} for more details). Following the same procedure of the Proof of Lemmas 1--3 in \cite{Li3}, we can prove \eqref{y2}--\eqref{y3} with suitable modification. Here we omit the details.

Now, we establish the estimate involving $\mathbf{S}_{t}(u_0)-u_0-t\mathbf{v}_0(u_0)$ which is crucial in proving Theorem \ref{the1.1}.
\begin{proposition}\label{pro1}
Assume that $||u_0||_{B^{\frac32}_{2,1}}\lesssim 1$. Under the assumptions of Theorem \ref{the1.1}, we have
\bal\label{et0}
||\mathbf{S}_{t}(u_0)-u_0-t\mathbf{v}_0(u_0)||_{B^{\frac32}_{2,1}}\leq Ct^{2}\mathbf{E}(u_0),
\end{align}
where we denote $\mathbf{v}_0(u_0):=\mathbf{P}(u_0)-u_0\pa_x u_0$ and
\bbal
\mathbf{E}(u_0)&:=1+||u_0||^2_{C^{0,1}}||u_0||_{B^{\frac52}_{2,1}}+||u_0||_{L^\infty}\Big(||u_0||_{B^{\frac52}_{2,1}}+\big(||u_0||_{L^\infty}
+||u_0||^2_{C^{0,1}}\big)||u_0||_{B^{\frac72}_{2,1}}\Big).
\end{align*}
\end{proposition}
{\bf Proof.}\quad For simplicity, we denote $u(t)=\mathbf{S}_t(u_0)$. Firstly, according to Lemma \ref{le5}, there exists a small time $T=T\big(||u_0||_{B^{\frac32}_{2,1}}\big)$ such that the solution $u(t)\in\mathcal{C}([0, T]; B^{\frac32}_{2,1})$, namely,
\bal\label{s}
||u(t)||_{L^\infty_TB^{\frac32}_{2,1}}\leq C||u_0||_{B^{\frac32}_{2,1}}\leq C.
\end{align}
Applying Lemma \ref{01} to Eq.\eqref{CH}, we have for all $t\in[0,T]$ and $\gamma\geq\frac32$
\bal\label{u2}
||u(t)||_{L^\infty_TB^\gamma_{2,1}}\leq ||u_0||_{B^{\gamma}_{2,1}}\exp\Big(C\int_0^T ||u||_{B^{\fr32}_{2,1}}\dd\tau\Big)\leq C||u_0||_{B^\gamma_{2,1}}.
\end{align}
Set
${\widetilde{u}}=\mathbf{S}_{t}(u_0)-u_0$, then we deduce from Eq.\eqref{CH} that
\begin{equation}\label{hy}
\begin{cases}
\pa_t{\widetilde{u}}+\mathbf{S}_{t}(u_0)\pa_x\mathbf{S}_{t}(u_0)=\mathbf{P}\big(\mathbf{S}_{t}(u_0)\big), \\
{\widetilde{u}}_0=0.
\end{cases}
\end{equation}
Then, we have
\bal\label{et1}
||\widetilde{u}(t)||_{L^\infty}&\leq\int^t_0||\pa_\tau \widetilde{u}||_{L^\infty} \dd\tau
\nonumber\\&\leq \int^t_0||\mathbf{S}_{t}(u_0)\pa_x\mathbf{S}_{t}(u_0)||_{L^\infty}\dd \tau+\int^t_0||\mathbf{P}\big(\mathbf{S}_{t}(u_0)||_{L^\infty}\dd \tau
\nonumber\\&\leq C\int^t_0||\mathbf{S}_{t}(u_0),u_0||^2_{C^{0,1}}\dd \tau\quad\text{by Lemma \ref{01}}\nonumber\\
&\leq Ct||u_0||^2_{C^{0,1}},
\end{align}
where we have used the estimate
$$||\mathbf{P}\big(\mathbf{S}_{t}(u_0)||_{L^\infty}\leq C||\mathbf{S}_{t}(u_0)||^2_{C^{0,1}}$$
from the fact that $(1-\pa^2_x)^{-1}f=G*f$ with $G(x)=\fr12e^{-|x|}$.

Note that $P(D)$ is a multiplier of degree $-1$, by Lemma \ref{le2}, we have for $\gamma\geq\fr32$
\bal\label{p}
||\mathbf{P}(u)||_{B^{\gamma}_{2,1}}&\leq C\Big\|u^2+\fr12u_x^2\Big\|_{B^{\gamma-1}_{2,1}}\nonumber\\
&\leq C\|u,u_x\|_{L^{\infty}}\|u,u_x\|_{B^{\gamma-1}_{2,1}}\nonumber\\
&\leq C||u||_{B^{\fr32}_{2,1}}||u||_{B^{\gamma}_{2,1}}.
\end{align}
By Lemma \ref{le2}, we obtain from \eqref{s} and \eqref{p} that
\bal\label{et2}
||u(t)-u_0||_{B^{\frac32}_{2,1}}
&\leq \int^t_0||\pa_\tau u||_{B^{\frac32}_{2,1}} \dd\tau
\nonumber\\&\leq \int^t_0||\mathbf{P}(u)||_{B^{\frac32}_{2,1}} \dd\tau+ \int^t_0||u \pa_xu||_{B^{\frac32}_{2,1}} \dd\tau
\nonumber\\&\leq Ct\Big(||u||^{2}_{B^{\frac32}_{2,1}}+||u||_{L^\infty}||u_x||_{B^{\frac32}_{2,1}}\Big)
\nonumber\\&\leq Ct\Big(||u_0||^{2}_{B^{\frac32}_{2,1}}+\big(||u_0||_{L^\infty}+||u_0||^2_{C^{0,1}}\big)||u_0||_{B^{\frac52}_{2,1}}\Big).
\end{align}

Similarly, we obtain from \eqref{s}, \eqref{u2} and \eqref{p} that
\bal\label{et3}
||u(t)-u_0||_{B^{\frac52}_{2,1}}
&\leq \int^t_0||\pa_\tau u||_{B^{\frac52}_{2,1}} \dd\tau
\nonumber\\&\leq \int^t_0||\mathbf{P}(u)||_{B^{\frac52}_{2,1}} \dd\tau+ \int^t_0||u \pa_xu||_{B^{\frac52}_{2,1}} \dd\tau
\nonumber\\&\leq Ct\Big(||u||_{B^{\frac32}_{2,1}}||u||_{B^{\frac52}_{2,1}}
+||u||_{L^\infty}||u||_{B^{\frac72}_{2,1}}\Big)
\nonumber\\&\leq Ct\Big(||u_0||_{B^{\frac32}_{2,1}}||u_0||_{B^{\frac52}_{2,1}}+\big(||u_0||_{L^\infty}+||u_0||^2_{C^{0,1}}\big)||u_0||_{B^{\frac72}_{2,1}}\Big).
\end{align}

Using Lemma \ref{le2} again, we obtain that
\bal\label{et4}
||u(t)-u_0-t\mathbf{v}_0(u_0)||_{B^{\frac32}_{2,1}}
&\leq \int^t_0||\pa_\tau u-\mathbf{v}_0(u_0)||_{B^{\frac32}_{2,1}} \dd\tau
\nonumber\\&\leq \int^t_0||\mathbf{P}(u)-\mathbf{P}(u_0)||_{B^{\frac32}_{2,1}} \dd\tau+\int^t_0||u\pa_xu-u_0\pa_xu_0||_{B^{\frac32}_{2,1}} \dd\tau
\nonumber\\&\lesssim \int^t_0||u(\tau)-u_0||_{B^{\frac32}_{2,1}} \dd\tau+\int^t_0||u(\tau)-u_0||_{L^\infty} ||u(\tau)||_{B^{\frac52}_{2,1}} \dd\tau
\nonumber\\&\quad \ + \int^t_0||u(\tau)-u_0||_{B^{\frac52}_{2,1}}  ||u_0||_{L^\infty}\dd \tau\nonumber\\
&\lesssim \int^t_0||u(\tau)-u_0||_{B^{\frac32}_{2,1}} \dd\tau+||u_0||_{B^{\frac52}_{2,1}}\int^t_0||u(\tau)-u_0||_{L^\infty}  \dd\tau
\nonumber\\&\quad \ + ||u_0||_{L^\infty}\int^t_0||u(\tau)-u_0||_{B^{\frac52}_{2,1}}  \dd \tau,
\end{align}
where we have used
\bbal
||\mathbf{P}(u)-\mathbf{P}(u_0)||_{B^{\fr32}_{2,1}}&\leq C\Big\|(u+u_0)(u-u_0)+\fr12\partial_x(u+u_0)\partial_x(u-u_0)\Big\|_{B^{\fr12}_{2,1}}\nonumber\\
&\leq C\|u-u_0\|_{B^{\fr32}_{2,1}}\|u+u_0\|_{B^{\fr32}_{2,1}}\nonumber\\
&\leq C\|u-u_0\|_{B^{\fr32}_{2,1}}.
\end{align*}
Plugging \eqref{et1}--\eqref{et2} into \eqref{et3} yields the desired result \eqref{et0}. Thus, we complete the proof of Proposition \ref{pro1}.

Also, the following estimate involving $\mathbf{S}_{t}(u_0)-u_0-t\mathbf{w}_0(u_0)$ is crucial in proving Theorem \ref{the1.2}.
\begin{proposition}\label{pro2}
Assume that $||u_0||_{B^{\frac32}_{2,1}}\lesssim 1$. Under the assumptions of Theorem \ref{the1.2}, we have
\bbal
&||\mathbf{S}_{t}(u_0)-u_0||_{L^\infty}\leq Ct||u_0||_{C^{0,1}},
\\&||\mathbf{S}_{t}(u_0)-u_0||_{B^{\frac32}_{2,1}}\leq Ct\big(||u_0||^{3}_{B^{\frac32}_{2,1}}
+||u_0||^2_{C^{0,1}}||u_0||_{B^{\frac52}_{2,1}}\big),
\\&||\mathbf{S}_{t}(u_0)-u_0||_{B^{\frac52}_{2,1}}\leq Ct\big(||u_0||^2_{B^{\frac32}_{2,1}}||u_0||_{B^{\frac52}_{2,1}}
+||u_0||^2_{C^{0,1}}||u_0||_{B^{\frac72}_{2,1}}\big),
\\&||\mathbf{S}_{t}(u_0)-u_0-t\mathbf{w}_0||_{B^{\frac32}_{2,1}}\leq Ct^{2}\mathbf{F}(u_0),
\end{align*}
where we denote $\mathbf{w}_0(u_0):=\mathbf{Q}(u_0)-u^2_0\pa_x u_0$ and
\bbal
\mathbf{F}(u_0):=1+||u_0||^2_{C^{0,1}}||u_0||_{B^{\frac52}_{2,1}}
+||u_0||^4_{C^{0,1}}||u_0||_{B^{\frac72}_{2,1}}.
\end{align*}
\end{proposition}
{\bf Proof.}\quad The proof follows the same manner from Proposition \ref{pro1}, we omit the details.
\section{Non-uniform continuous dependence}

In this section, we prove Theorems \ref{the1.1} and \ref{the1.2} by using Propositions \ref{pro1} and \ref{pro2}, respectively.

{\bf Proof of Theorem \ref{the1.1}}\quad
We set $u^n_0=f_n+g_n$ and compare the solution $\mathbf{S}_{t}(u^n_0)$ and $\mathbf{S}_{t}(f_n)$. We obviously have
\bbal
||u^n_0-f_n||_{B^{\frac32}_{2,1}}=||g_n||_{B^{\frac32}_{2,1}}\leq C2^{-n},
\end{align*}
which means that
\bbal
\lim_{n\to\infty}||u^n_0-f_n||_{B^{\frac32}_{2,1}}=0.
\end{align*}
From Lemma \ref{le4}, one has
\bbal
&||u^n_0,f_n||_{B^{\sigma}_{2,1}}\leq C2^{(\sigma-\fr32)n}\quad \text{for}\quad \sigma\geq\fr32,\\
&\|u^n_0,f_n\|_{L^\infty}\leq C2^{-n}\quad\text{and}\quad\|\pa_xu^n_0,\pa_xf_n\|_{L^\infty}\leq C2^{-\fr{n}{2}},
\end{align*}
which implies
\bbal
\mathbf{E}(u^n_0)+\mathbf{E}(f_n)\leq C.
\end{align*}
Notice that
\bbal
&\mathbf{S}_{t}(u^n_0)=\underbrace{\mathbf{S}_{t}(u^n_0)-u^n_0-t\vv_0(u_0^n)}_{=~\mathbf{I}_1(\vv_0)}+f_n+g_n+t\big(\mathbf{P}(u^n_0)-u^n_{0}\pa_xu^n_{0}\big)\nonumber\\
&\mathbf{S}_{t}(f_n)=\underbrace{\mathbf{S}_{t}(f_n)-f_n-t\vv_0(f_n)}_{=~\mathbf{I}_2(\vv_0)}+f_n+t\big(\mathbf{P}(f_n)-f_n\pa_x f_n\big)\quad\text{and}\nonumber\\
&u^n_{0}\pa_xu^n_{0}-f_n\pa_xf_n=g_n\pa_xf_n+u^n_{0}\pa_xg_n,
\end{align*}
using the triangle inequality and Proposition \ref{pro1}, we deduce that
\bal\label{yyh}
\quad \ &\big\|\mathbf{S}_{t}(u^n_0)-\mathbf{S}_{t}(f_n)\big\|_{B^{\frac32}_{2,1}}\nonumber\\
=&~\big\|\mathbf{I}_1(\vv_0)-\mathbf{I}_2(\vv_0)+g_n-t\big(g_n\pa_xf_n+u^n_{0}\pa_xg_n-\mathbf{P}(u^n_0)+\mathbf{P}(f_n)\big)\big\|_{B^{\frac32}_{2,1}}\nonumber\\
\geq&~t\big\|g_n\pa_xf_n+u^n_{0}\pa_xg_n-\big(\mathbf{P}(u^n_0)-\mathbf{P}(f_n)\big)\big\|_{B^{\frac32}_{2,1}}-\big\|\mathbf{I}_1(\vv_0),\mathbf{I}_2(\vv_0),g_n\big\|_{B^{\frac32}_{2,1}}\nonumber\\
\geq&~ t\big\|g_n\pa_xf_n\big\|_{B^{\frac32}_{2,1}}-t\big\|u^n_{0}\pa_xg_n,\;\mathbf{P}(u^n_0)-\mathbf{P}(f_n)\big\|_{B^{\frac32}_{2,1}}-Ct^{2}-C2^{-n}\nonumber\\
\geq&~ t\big\|g_n\pa_xf_n\big\|_{B^{\frac32}_{2,1}}-Ct2^{-n}-Ct^{2}-C2^{-n},
\end{align}
where we have performed the following easy computations
\bbal
&\big\|u^n_{0}\pa_xg_n\big\|_{B^{\frac32}_{2,1}}\leq C\big\|u^n_0\big\|_{B^{\frac32}_{2,1}}\big\|g_n\big\|_{B^{\frac52}_{2,1}}\leq C2^{-n},\\
&\big\|\mathbf{P}(u^n_0)-\mathbf{P}(f_n)\big\|_{B^{\frac32}_{2,1}}\leq C\big\|g_n\big\|_{B^{\frac32}_{2,1}}\big\|u^n_0+f_n\big\|_{B^{\frac32}_{2,1}}\leq C2^{-n}.
\end{align*}
Combining the fact from Lemma \ref{le4}
\begin{eqnarray*}
      \liminf_{n\rightarrow \infty} \big\|g_n\pa_xf_n\big\|_{B^{\frac32}_{2,1}}\gtrsim M_1,
        \end{eqnarray*}
then we deduce from \eqref{yyh} that
\bbal
\liminf_{n\rightarrow \infty}\big\|\mathbf{S}_t(f_n+g_n)-\mathbf{S}_t(f_n)\big\|_{B^{\frac32}_{2,1}}\gtrsim t\quad\text{for} \ t \ \text{small enough}.
\end{align*}
This completes the proof of Theorem \ref{the1.1}.

{\bf Proof of Theorem \ref{the1.2}}\quad
We set $u^n_0=f_n+h_n$ and compare the solution $\mathbf{S}_{t}(u^n_0)$ and $\mathbf{S}_{t}(f_n)$. Obviously, we have
\bbal
\lim_{n\to\infty}\|u^n_0-f_n\|_{B^{\frac32}_{2,1}}=\lim_{n\to\infty}\|h_n\|_{B^{\frac32}_{2,1}}=0.
\end{align*}
Lemma \ref{le4} tells us that
\bbal
\|u^n_0,f_n\|_{C^{0,1}}\leq C2^{-\frac{n}{2}} \quad  \mathrm{and}\quad
\|u^n_0,f_n\|_{B^{\sigma}_{2,1}}\leq C2^{(\sigma-\frac32)n}\quad  \mathrm{for} \quad \sigma\geq \frac32,
\end{align*}
which implies
\bbal
\mathbf{F}(u^n_0)+\mathbf{F}(f_n)\leq C.
\end{align*}
Using the triangle inequality and Proposition \ref{pro2}, we deduce that
\bal\label{yyh2}
\quad \ &\|\mathbf{S}_{t}(u^n_0)-\mathbf{S}_{t}(f_n)\|_{B^{\frac32}_{2,1}}\nonumber\\
=&~\big\|\mathbf{I}_1(\mathbf{w}_0)-\mathbf{I}_2(\mathbf{w}_0)+g_n-t\big((u^n_{0})^2\pa_xu^n_{0}-f^2_n\pa_xf_n-\mathbf{Q}(u^n_0)+\mathbf{Q}(f_n)\big)\big\|_{B^{\frac32}_{2,1}}\nonumber\\
\geq&~t\big\|h^2_n\pa_xf_n\big\|_{B^{\frac32}_{2,1}}-t\|2f_nh_n\pa_xf_n+(u^n_{0})^2\pa_xh_n+\mathbf{Q}(f_n)
-\mathbf{Q}(u^n_0)\big\|_{B^{\frac32}_{2,1}}\nonumber\\
&-\big\|\mathbf{I}_1(\mathbf{w}_0),\mathbf{I}_2(\mathbf{w}_0),g_n\big\|_{B^{\frac32}_{2,1}}\nonumber\\
\geq&~ t\big\|h^2_n\pa_xf_n\big\|_{B^{\frac32}_{2,1}}-Ct2^{-\frac{n}{2}}-Ct^{2},
\end{align}
where we have used that
$$
(u^n_{0})^2\pa_xu^n_{0}-f^2_n\pa_xf_n=h^2_n\pa_xf_n+2f_nh_n\pa_xf_n+(u^n_{0})^2\pa_xh_n
$$
and
\bbal
&\big|\big|(u^n_{0})^2\pa_xh_n\big|\big|_{B^{\frac32}_{2,1}}\leq C||u^n_0||^2_{B^{\frac32}_{2,1}}||h_n||_{B^{\frac52}_{2,1}}\leq C2^{-\frac{n}{2}},
\\&\big|\big|f_nh_n\pa_xf_n\big|\big|_{B^{\frac32}_{2,1}}\lesssim ||f_n||_{L^\infty}||h_n||_{L^\infty}||f_n||_{B^{\frac52}_{2,1}}+||\pa_xf_n||_{L^\infty}||h_n||_{B^{\frac32}_{2,1}}||f_n||_{B^{\frac32}_{2,1}}\leq C2^{-n}.
\end{align*}
Combining the fact from Lemma \ref{le4}
\begin{eqnarray*}
      \liminf_{n\rightarrow \infty} \big\|h^2_n\pa_xf_n\big\|_{B^{\frac32}_{2,1}}\gtrsim M_2,
        \end{eqnarray*}
then we deduce from \eqref{yyh2} that
\bbal
\liminf_{n\rightarrow \infty}||\mathbf{S}_t(f_n+h_n)-\mathbf{S}_t(f_n)||_{B^{\frac32}_{2,1}}\gtrsim t\quad\text{for} \ t \ \text{small enough}.
\end{align*}
This completes the proof of Theorem \ref{the1.2}.

\vspace*{1em}
\noindent\textbf{Acknowledgements.}  J. Li is supported by the National Natural Science Foundation of China (Grant No.11801090). Y. Yu is supported by the Natural Science Foundation of Anhui Province (No.1908085QA05). W. Zhu is partially supported by the National Natural Science Foundation of China (Grant No.11901092) and Natural Science Foundation of Guangdong Province (No.2017A030310634).

\end{document}